\documentclass[11pt,leqno]{article} 
\usepackage{graphics}
\newtheorem{thm}{Theorem}[section]
\newtheorem{lma}{Lemma}[section]

\newcommand{\beqa}{\begin{eqnarray}}
\newcommand{\eeqa}{\end{eqnarray}}

\newcommand{\pf}{\noindent {\bf Proof:} $\s$ }
\newcommand{\epf}{ \hfill$\diamondsuit$ \medskip}

\newcommand{\beq}{\begin{equation}}
\newcommand{\eeq}{\end{equation}}
\newcommand{\lbl}{\label}
\newcommand{\s}{\; \;}

\newcommand{\la}{\lambda}

\newcommand{\ra}{\rightarrow}

\title{On the interaction of species capable of explosive growth}

\author{
Philip Korman  \\ 
Department of Mathematical Sciences \\ 
University of Cincinnati \\ 
Cincinnati Ohio 45221-0025 \\
}

\date{}

\begin{document}

\maketitle
\begin{abstract} 
In the classical   Lotka-Volterra population models, the interacting species affect each other's growth rate. We propose an alternative model, in which the species affect each other through the limitation coefficients, rather  then through the growth rates. This appears to be more realistic: the presence of foxes is not likely to diminish the fertility of rabbits, but will contribute to limiting rabbit's population.  Both the cases of predation and of competition are considered, as well as competition in case of periodic coefficients. Our model becomes linear when one switches to the reciprocals of the variables.
In another direction we use a similar idea to derive a multiplicity result for a class of periodic equations.
 \end{abstract}

\begin{flushleft}
\end{flushleft}

\begin{flushleft}
\end{flushleft}

\section{Introduction}
\setcounter{equation}{0}
\setcounter{thm}{0}
\setcounter{lma}{0}

One way of solving the logistic population equation (here $x=x(t)$)
\beq
\lbl{l1}
x'=ax-bx^2
\eeq
is to divide this equation by $x^2$, and obtain a linear equation for $u=\frac{1}{x}$. Here $a>0$ is the growth rate, and $b>0$ is the limitation (or self-limitation) coefficient, both given numbers. We wish to explore the interactions of two species with populations $x=x(t)$ and $y=y(t)$ for which the substitution $u=\frac{1}{x}$ and $v=\frac{1}{y}$ leads to a linear system. The model we consider is
\beqa
\lbl{l2}
& x'=ax+x^2 \left(\frac{b}{y}+e \right) \\ \nonumber
& y'=dy+y^2 \left(\frac{c}{x}+f \right) \,,
\eeqa
with constants $a$,$b$,$c$,$d$,$e$ and $f$. Dividing the first  equation by $x^2$, the second one by $y^2$, and setting $u=\frac{1}{x}$ and $v=\frac{1}{y}$, gives  a linear system
\beqa 
\lbl{12b}
& -u'=au+bv+e \\ \nonumber
& -v'=cu+dv+f \,.
\eeqa
The signs of the coefficients determine the type of interaction, which will include both predator-prey and competing species cases.
\medskip

Let us compare  (\ref{l2}) with the classical Lotka-Volterra predator-prey model 
\beqa 
\lbl{l2a}
& x'=x \left(a-b\, y \right) \\ \nonumber
& y'=y \left(-c+d \, x \right) \,,
\eeqa
where the constants $a$,$b$,$c$,$d$ are  positive. In (\ref{l2a}) the species affect each other through the growth rate: the prey, with the number given by $x(t)$, improves the growth rate of the predator, with the number $y(t)$, while the predator decreases the growth rate of the prey. In the model (\ref{l2}) the species affect each other through their limitation coefficients. This appears to be more realistic: the presence of foxes is not likely to decrease the fertility of rabbits (new  rabbits will be born at the same rate), but will place a limitation on the  growth of rabbit population. 
\medskip

Similarly to the Lotka-Volterra  model, the proposed  model (\ref{l2}) predicts oscillatory behavior for predator-prey interaction, and either stable coexistence or competitive exclusion for competing species. Unlike  the Lotka-Volterra  model, it is possible that the population number of one of the species  goes to infinity in finite time, while the number of the other species remains finite and positive. Explosive growth of populations occurs often in nature. Notice that our analysis leads to some non-standard questions about linear systems. For example, if a solution of (\ref{12b}) starts in the first quadrant of the $xy$-plane, will it stay in the first quadrant for all $t$?
\medskip

Using the Floquet theory, we analyze a case of predator-prey interaction with periodic coefficients, and give a condition for the existence of a  limit cycle.
\medskip

In another direction we use the same transformation $ u= \frac{1}{x}$  to derive a multiplicity result for a class of periodic equations
\[
x'(t)=f(t,x(t)) \,, \s\s \mbox{with $f(t+p,x)=f(t,x)$} \,.
\]

\section{Explosive predator-prey model}
\setcounter{equation}{0}
\setcounter{thm}{0}
\setcounter{lma}{0}

Consider the model
\beqa
\lbl{l3}
& x'=x^2 \left(\frac{b}{y}-1 \right) \\ \nonumber
& y'=-y^2 \left(\frac{d}{x}-1 \right) \,.
\eeqa
Here $x(t)$ gives the number of prey, and $y(t)$  the number of predator. If $y(t)$ is small, the prey  grows explosively (with $x'$ behaving like $\alpha x^2$, $\alpha >0$). If the number of predators $y(t)$ is large, then $x'(t)<0$ and $x(t)$ decreases. The number of  predators $y(t)$ decreases when $x(t)$ small, and  grows explosively for $x(t)$ large. This model corresponds to (\ref{l2}), with $a=d=0$. The coefficients $e$ and $f$ have been scaled out.
\medskip

The system (\ref{l3}) has a rest point $(d,b)$. Letting $X=\frac{\sqrt{d/b}}{x}$ and $Y=\frac{1}{y}$ transforms (\ref{l3}) into  a perturbed harmonic oscillator
\beqa 
\lbl{l4}
& X'=-\sqrt{bd} \, Y+\sqrt{d/b} \\ \nonumber
& Y'=\sqrt{bd}\, X-1 \,.
\eeqa
Setting $X(t)=\xi(t)+\frac{1}{\sqrt{bd}}$, $Y(t)=\eta(t)+\frac{1}{b}$ leads to a harmonic oscillator
\beqa \nonumber
& \xi '=-\sqrt{bd} \, \eta \\ \nonumber
& \eta '=\sqrt{bd}\, \xi \,,
\eeqa
so that the solution of (\ref{l4}) is
\beqa 
\lbl{14.1}
& X(t)=\frac{1}{\sqrt{bd}}+c_1 \cos \sqrt{bd} \, t-c_2 \sin \sqrt{bd} \, t \\ \nonumber
& Y(t)=\frac{1}{b}+c_1 \sin \sqrt{bd} \, t+c_2 \cos  \sqrt{bd} \, t \,,
\eeqa
which is just a rotation of the point $(X(0),Y(0))$ around the point $(\frac{1}{\sqrt{bd}},\frac{1}{b})$, the rest point of (\ref{l4}),  on the circle of radius $\sqrt{c_1^2+c_2^2}$. 
The solution of (\ref{l3}) is then
\beqa 
\lbl{l6}
& x(t)=\frac{\sqrt{d/b}}{\frac{1}{\sqrt{bd}}+c_1 \cos \sqrt{bd} \, t-c_2 \sin \sqrt{bd} \, t} \\ \nonumber
& y(t)=\frac{1}{\frac{1}{b}+c_1 \sin \sqrt{bd} \, t+c_2 \cos  \sqrt{bd} \, t} \,.
\eeqa
The constants $c_1$ and $c_2$ are determined from the initial values $(x(0),y(0))$: 
\beq
\lbl{l7}
c_1=\frac{\sqrt{d/b}}{x(0)}-\frac{1}{\sqrt{bd}} \,, \s \mbox{ and} \s c_2=\frac{1}{y(0)}-\frac{1}{b} \,.
\eeq
It is now clear that the rest point  $(d,b)$ is a center for (\ref{l3}), and we can  give a complete description of the behavior of positive solutions.

\begin{thm}
Given the initial point $(x(0),y(0))$, calculate $c_1$ and $c_2$ by (\ref{l7}), and $R=\sqrt{c_1^2+c_2^2}$. If the circle $C$ of radius $R$ around  the point $(\frac{1}{\sqrt{bd}},\frac{1}{b})$ lies completely inside the first quadrant of the $(X,Y)$ plane, then the corresponding solution $(x(t),y(t))$ of (\ref{l3}) is a closed curve around the rest point $(d,b)$, given by (\ref{l6}). Moreover, the period of all these closed curves is the same, and $x(t) > \frac{d}{2}$, $y(t)> \frac{b}{2}$ for all $t$. Assume now that this circle $C$, traveled  counterclockwise beginning with the point $(X(0),Y(0))=(\frac{\sqrt{d/b}}{x(0)},\frac{1}{y(0)})$, hits one of the axes  of the $(X,Y)$ plane. If it hits the $Y$-axis first, then there is a time $T>0$ so that $\lim _{t \ra T} x(t)=\infty$, while   $\lim _{t \ra T} y(t)$ is finite and positive. If $C$ hits the $X$-axis first, then there is a time $T>0$ so that $\lim _{t \ra T} y(t)=\infty$, while   $\lim _{t \ra T} x(t)$ is finite and positive. 
\end{thm}

\pf
In view of the discussion above, it remains to prove the lower bounds for the periodic solutions in the first part of the theorem. From (\ref{14.1}) one sees that the positivity of $X(t)$ and $Y(t)$ implies that $X(t)<\frac{2}{\sqrt{bd}}$ and $Y<\frac{2}{b}$, from which one gets the lower bounds on $x(t)$ and $y(t)$.
\epf

\noindent
{\bf Example} Using {\em Mathematica}, we computed four periodic solutions for   the system (\ref{l3}), with $b=3$ and $d=2$, surrounding the rest point at $(2,3)$, see Figure \ref{fig:2}.

\begin{figure}
\scalebox{0.65}{\includegraphics{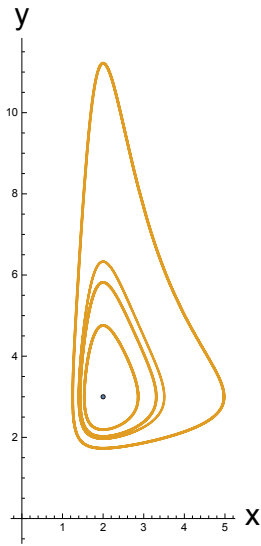}}
\caption{  Periodic solutions  for the system (\ref{l3})}
\lbl{fig:2}
\end{figure}

\section{Explosive competing species model}
\setcounter{equation}{0}
\setcounter{thm}{0}
\setcounter{lma}{0}

Consider the model
\beqa
\lbl{20}
& x'=a\, x+x^2 \left(\frac{b}{y}-1 \right),  \s x(0)>0 \\ \nonumber
& y'=d \, y+y^2 \left(\frac{c}{x}-1 \right), \s y(0)>0 \,,
\eeqa
with positive constants $a$,$b$,$c$ and $d$. Each species  grows explosively, if the number of the other one is small, while if the competitor's  number is large, the growth is logistic-like. Clearly, the interaction is competitive in nature.
\medskip

We begin with a simple observation: if $x(0)>0$ and $y(0)>0$, then $x(t)>0$ and $y(t)>0$ for all $t>0$. Indeed, writing the first equation in the form $x'=A(t)x$, with $A(t) \equiv a+x(t) \left(\frac{b}{y(t)}-1 \right)$, and integrating, obtain $x(t)=x(0)e^{\int_0^t A(s) \, ds}>0$. Similarly, $y(t)>0$ for all $t>0$. Hence, we can limit our study of (\ref{20}) to the first quadrant of the $(x,y)$ plane.
\medskip

Setting $X=\frac{1}{x}$ and $Y=\frac{1}{y}$ produces a linear system
\beqa 
\lbl{21}
& X'=-aX-bY+1 \,, \s X(0)=\frac{1}{x(0)}>0 \\ \nonumber
& Y'=-cX-dY+1 \,,  \s Y(0)=\frac{1}{y(0)}>0 \,,
\eeqa
with a unique rest point $(X_0,Y_0)$ given by
\beq
\lbl{22}
X_0=\frac{d-b}{ad-bc} \,, \s Y_0=\frac{a-c}{ad-bc} \,.
\eeq
Since $x(t)>0$ and $y(t)>0$ for all $t>0$, we may restrict the system (\ref{21}) to the first quadrant of the $(X,Y)$ plane.
The  rest point $(X_0,Y_0)$ lies in the first quadrant if either 
\beq
\lbl{23}
d>b \; \mbox{ and} \; a>c \; \mbox{ (and then $ad>bc$)} \,, 
\eeq
or 
\beq
\lbl{24}
d<b \; \mbox{ and} \; a<c \; \mbox{ (and then $ad<bc$)} \,. 
\eeq
Letting $\xi=X-X_0$ and $\eta=Y-Y_0$, we translate the rest point to the origin, obtaining the system
\beqa 
\lbl{25}
& \xi '=-a\xi-b\eta \\ \nonumber
& \eta '=-c\xi-d\eta \,,
\eeqa
with the matrix $A=
\left[
\begin{array}{rr}
-a & -b \\
-c & -d
\end{array}
\right]$. The eigenvalues of $A$ are 
\[
\lambda _{1,2}=\frac{1}{2} \left(-a-d \pm \sqrt{a^2-2
   a d+4 b c+d^2}\right) \,.
\]
The corresponding (column) eigenvectors are 
\[
\xi  _{1,2}= \left( -\frac{-a+d \pm\sqrt{a^2-2 a d+4 b
   c+d^2}}{2 c},1 \right)^T \,.
\]
In case (\ref{23}) holds,  both eigenvalues are negative, and the  rest point $(X_0,Y_0)$ is a stable node, while in case (\ref{24}) holds, one eigenvalue is negative  and the other one is positive, so that $(X_0,Y_0)$ is a saddle.

\begin{thm}
(i) Assume that the condition  (\ref{24}) holds. Then one of the species (depending on the initial conditions) grows explosively. Namely, for any solution of (\ref{20}) there is a time $T>0$ so that $\lim _{t \ra T} x(t)=\infty$, while   $\lim _{t \ra T} y(t)$ is finite and positive, or the other way around. 
\medskip

\noindent
(ii) Assume that the condition  (\ref{23}) holds. 
If $x(t)$ and $y(t)$ remain finite for all $t>0$ then $\lim _{t \ra \infty} x(t)=\frac{1}{X_0}$ and $\lim _{t \ra \infty} y(t)=\frac{1}{Y_0}$.
\end{thm}

\pf
The general solution of (\ref{21}) is
\beq
\lbl{28a}
\left(X(t),Y(t) \right)^T=\left(X_0,Y_0 \right)^T+c_1e^{\la _1t}  \xi _1+c_2e^{\la _2t}  \xi _2 \,.
\eeq
(i) In case (\ref{24}) holds,  the eigenvalues of $A$ are of opposite sign say $\la _1<0<\la _2$. The term $c_1e^{\la _1t}  \xi _1$ is negligible in the long run. The eigenvector $\xi _2$ corresponding to the positive eigenvalue (``plus" in front of the square root) has one component positive, and the other one is negative. It follows that all of the solutions of (\ref{21}) eventually move either northwest or southeast of the rest point $(X_0,Y_0)$ intersecting either the $X$ or the $Y$ axis.  
\medskip

\noindent
(ii) In case (\ref{23}) holds, the general solution of (\ref{21}) is
given by (\ref{28a}), 
with negative $\la _1$ and $\la _2$. It follows that the point $\left(X(t),Y(t) \right)$  tends to the point $\left(X_0>0,Y_0>0 \right)$ as $t \ra \infty$.
If the  point $\left(X(t),Y(t) \right)$ stays in the first quadrant, then $x(t)$ and $y(t)$ are defined for all $t$, otherwise one of the species becomes infinite in finite time.
\epf

\noindent
{\bf Remark}  $\;$ In case (\ref{23}) holds, the solution of  (\ref{21}) connects the points $\left(X_0,Y_0 \right)$ and $\left(X(0),Y(0 )\right)$ in the first quadrant.  While it is rare for the solution   $\left(X(t),Y(t) \right)$ to exit the first quadrant, this may indeed happen if the points $\left(X_0,Y_0 \right)$ and $\left(X(0),Y(0 )\right)$ lie near one of the axes. We used {\em Mathematica} to solve (\ref{21}) with $a=4$, $b=1$, $c=1$, $d=5$, $X(0)=5$, $Y(0)=0.1$. Here $X_0=\frac{4}{19}>0$ and $Y_0=\frac{3}{19}>0$. The graph of the solution in Figure \ref{fig:1} shows that $Y(t)$ becomes zero at some $T$, which corresponds to $\lim _{t \ra T}  y(t)=\infty$.
\begin{figure}[h]
\begin{center}
\scalebox{0.99}{\includegraphics{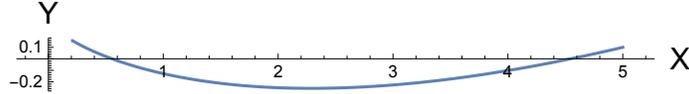}}
\caption{  A  solution of the system (\ref{21}) exiting the first quadrant of the $XY$-plane (the motion is from right to left)}
\lbl{fig:1}
\end{center}
\end{figure}
 
\section{Explosive predator-prey model  with periodic coefficients}
\setcounter{equation}{0}
\setcounter{thm}{0}
\setcounter{lma}{0}

We now consider a periodic perturbation of the explosive predator-prey model
\beqa
\lbl{30}
& x'=x^2 \left(\frac{b+\beta (t)}{y}-1 \right) \\ \nonumber
& y'=-y^2 \left(\frac{d+\delta (t)}{x}-1 \right) \,,
\eeqa
with small continuous functions $\beta (t)$ and $\delta (t)$ of period $p$, so that $\beta (t+p)=\beta (t)$ and $\delta (t+p)=\delta (t)$ for all $t$. (We make no assumptions on the sign of $\beta (t)$ and $\delta (t)$.) The linear system for $X=\frac{1}{x}$ and $Y=\frac{1}{y}$
\beqa
\lbl{31}
& X'= -\left(b+\beta (t) \right)Y+1 \\ \nonumber
& Y'= \left(d+\delta (t) \right)X-1 
\eeqa
has $p$-periodic coefficients. Let $F(t)$ be the normalized fundamental solution matrix (with $F(0)=I$, the identity matrix) of the corresponding homogeneous system 
\beqa
\lbl{32}
& X'= -\left(b+\beta (t) \right)Y \\ \nonumber
& Y'= \left(d+\delta (t) \right)X \,. 
\eeqa
For small  $\beta (t)$ and $\delta (t)$, $F(t)$ is close for $t \in [0,p]$ to the normalized fundamental solution matrix $F_0(t)=
\left[
\begin{array}{rr}
\cos \sqrt{bd} \, t & - \sqrt{\frac{b}{d}}\sin \sqrt{bd}  \, t \\
\sqrt{\frac{d}{b}} \sin \sqrt{bd}  \, t & \cos \sqrt{bd} \, t
\end{array}
\right]
$ of the unperturbed system
\beqa
\lbl{33}
& X'= -bY \\ \nonumber
& Y'= dX \,. 
\eeqa
By the continuous dependence of eigenvalues on the coefficients of the matrix, the Floquet multipliers of (\ref{32}), i.e., the eigenvalues of  $F(p)$ are close to the eigenvalues $\rho _1$ and $\rho _2$ of $F_0(p)$. Clearly,
\beqa
\lbl{34}
& \rho _1  \rho _2=1=\det F_0(p)\\ \nonumber
& \rho _1 + \rho _2= 2 \cos \sqrt{bd} \, p= {\rm trace} \, F_0(p)\,. 
\eeqa

\begin{thm}
Assume that $\sqrt{bd} \, p \ne 2 \pi m$, for any integer $m$. Then the system (\ref{30}) has a unique positive $p$-periodic solution $(x_p(t),y_p(t))$ for sufficiently small  $\beta (t)$ and $\delta (t)$. 
\end{thm}

\pf
Observe that $\rho _i \ne 1$, for $i=1,2$. Indeed, if $\rho_1=1$, then from the first line in (\ref{34}) $\rho_2=1$, giving a contradiction in the second line in (\ref{34}), because $\cos \sqrt{bd} \, p \ne 1$. Since $\beta (t)$ and $\delta (t)$ are small, the Floquet multipliers of the homogeneous problem (\ref{32}) are different from one, so that   (\ref{32}) has no $p$-periodic solution, and then by a standard result the non-homogeneous system  (\ref{31}) (and hence the original system (\ref{30})) has a unique $p$-periodic solution   $(X_p(t),Y_p(t))$. It remains to show that $X_p(t)>0$ and $Y_p(t)>0$ for all $t$.
\medskip

We derive  next an a priori bound on $X_p(t)$ and $Y_p(t)$, uniform in $\beta (t)$ and $\delta (t)$, provided that $|\beta (t)|+|\delta (t)| \leq c_0$, for some constant $c_0$.  Indeed, integrating both equations in (\ref{31}) over $(0,t)$, with $t \in (0,p)$, taking absolute values and then adding the corresponding inequalities, obtain
\beq
\lbl{35.1}
|X_p(t)|+|Y_p(t)| \leq a_1\int_0^t \left(|X_p(s)|+|Y_p(s)| \right) \, ds+a_2 \,,
\eeq
for some positive constants $a_1$ and $a_2$. The desired bound over $(0,p)$ follows by the Bellman-Gronwall lemma, see e.g., \cite{K}. 
\medskip

We claim that $X_p(t)>0$ and $Y_p(t)>0$ for all $t$.  Setting $X_p(t)=\xi (t)+\frac{1}{d}$ and $Y_p(t)=\eta (t)+\frac{1}{b}$ in  (\ref{31}) obtain
\beqa 
\lbl{36}
& \xi '=-b \, \eta -\beta (t) Y_p(t)\\ \nonumber
& \eta '=d \xi +\delta (t) X_p(t)\,.
\eeqa
Express
\[
\left[
\begin{array}{r}
\xi(t) \\
\eta  (t)
\end{array}\right]=F_0(t)
\left[
\begin{array}{r}
c_1 \\
c_2
\end{array}\right]+ F_0(t) \int_0^t F_0^{-1}(s) f(s) \, ds \,,
\]
with some constants $c_1$ and $c_2$, and $ f(t)=
\left[
\begin{array}{r}
-\beta (t) Y_p(t) \\
\delta (t) X_p(t)
\end{array}\right]$. Since the vector $\left[
\begin{array}{r}
\xi(t) \\
\eta  (t)
\end{array} \right]$ has period $p$, and the fundamental solution matrix $F_0(t)$ has period $\frac{2 \pi m}{\sqrt{bd}} \ne p$, it follows that $c_1=c_2=0$.
The vector 
$ f(t)$ is small by our assumptions, and the a priori estimate (\ref{35.1}). Both matrices $F_0(t)$ and $F_0^{-1}(s)$ have bounded entries. Then  the vector $\left[
\begin{array}{r}
\xi(t) \\
\eta  (t)
\end{array} \right]$ is small, so that the trajectory $(X_p(t),Y_p(t))$ remains near the point $\left( \frac{1}{d},\frac{1}{b} \right)$, and hence it stays in the first quadrant for all $t$.
\epf

\section{Multiplicity of solutions for a class of periodic equations}
\setcounter{equation}{0}
\setcounter{thm}{0}
\setcounter{lma}{0}

The transformation $u(t)=\frac{1}{x(t)}$ of the preceding sections turns out to be useful for a class of first order equations with periodic coefficients. V.A. Pliss \cite{P} considered what he called the {\em Abel equation}:
\beq
\lbl{50}
x'(t)=a_0(t)x^3+a_1(t)x^2+a_2(t)x+a_3(t) \,.
\eeq
Assuming that the given functions $a_i(t)$, $0 \leq i \leq 3$, are of period $p$, and $a_0(t)$ is either positive or negative for all $t$,  he proved that the equation (\ref{50}) has at most three $p$-periodic solutions. The proof involved a clever combination of the equations that the inverses of solutions satisfy.
\medskip

What if one changes the $a_0(t)x^3$ term to $a_0(t)x^{2n+1}$? In case it is $a_0(t)x^5$, the method of  V.A. Pliss \cite{P} still gives the same result with a little extra effort. For higher powers things get more involved, and in fact existence of at most three $p$-periodic solutions was proved by another elegant method in A.A. Panov \cite{P1}. It turns out that the following more general result was already known.

\begin{thm}\lbl{thm:5.1}
For the equation 
\beq
\lbl{51}
x'(t)=f(t,x) 
\eeq
assume that   the function $f(x,t)$ is continuous and has three continuous derivatives in $x$, and also for some $p>0$ and  all real $t$ and $x$ one has
\beq
\lbl{51a}
f(t+p,x)=f(t,x) \,,
\eeq
\beq
\lbl{52}
f_{xxx}(t,x)>0 \s\s \mbox{(or the opposite inequality holds)} \,.
\eeq
Then the equation (\ref{51}) has at most three $p$-periodic solutions.
\end{thm}

This theorem follows from a more general result of A.  Sandqvist and K.M. Andersen \cite{S}. They considered the equation (\ref{51}) on the interval $(0,p)$ and called a solution to be {\em closed} if $x(p)=x(0)$. Assuming the condition (\ref{52}) holds, they showed that the problem (\ref{51}) has at most three closed solutions, which implies the Theorem \ref{thm:5.1}.
\medskip

A simpler proof of the Theorem \ref{thm:5.1} was found in P. Korman and T. Ouyang \cite{KO}. We now simplify the presentation in that paper. The proof will follow from the following three simple lemmas.

\begin{lma}\lbl{lma:5.1}
Assume the condition (\ref{52}) holds and $f(t,0)=0$ for all $t \in R$. Then for all $t \in R$ and $x>0$ one has
\[
Q(t,x) \equiv 2f(t,x)-2x f_x(t,x)+x^2 f_{xx}(t,x) >0 \,.
\]
\end{lma}

\pf
Calculate $Q(t,0)=0$ and $Q_x(t,x)=x^2 f_{xxx}(t,x) >0$.
\epf

\begin{lma}\lbl{lma:5.2}
For the problem
\beq
\lbl{53}
y'(t)=g(t,y)
\eeq
assume that for some $p>0$ and all $t \in R$ and $y>0$ one has
\[
g(t+p,y)=g(t,y) \,,
\]
\[
g_{yy}(t,y)>0 \s\s \mbox{(or the opposite inequality holds)} \,.
\]
Then the equation (\ref{53}) has at most two positive $p$-periodic solutions.
\end{lma}

The proof is standard, and it can be found in e.g., P. Korman \cite{K}, p. $245$. The next lemma is crucial.

\begin{lma}\lbl{lma:5.3}
For the problem (\ref{51}) assume that $f(t,0)=0$ for all $t \in R$,   and the conditions (\ref{51a}),(\ref{52}) hold for all $t \in R$ and $x>0$. Then the equation (\ref{51}) has at most two positive $p$-periodic solutions.
\end{lma}

\pf
Set $x=\frac{1}{y}$ in  (\ref{51}). Then
\beq
\lbl{54}
-y'=y^2 f(t,\frac{1}{y}) \equiv g(t,y).
\eeq
By Lemma \ref{lma:5.1} for any $y>0$
\[
g_{yy}=2 f(t,\frac{1}{y})-\frac{2}{y}f_x(t,\frac{1}{y})+\frac{1}{y^2}f_{xx}(t,\frac{1}{y})=2f(t,x)-2x f_x(t,x)+x^2 f_{xx}(t,x) >0 \,.
\]
By Lemma \ref{lma:5.2} the equation  (\ref{54}) has at most two positive $p$-periodic solutions, and the same is true for (\ref{51}).
\epf

Turning to the proof of the Theorem \ref{thm:5.1}, observe that different solutions of (\ref{51}) do not intersect by the uniqueness theorem. If the equation  (\ref{51})  has four $p$-periodic solutions, let $\xi (t)$ be the smallest one. Then $z(t)=x(t)-\xi (t)$ satisfies 
\beq
\lbl{55}
z'=f(t,z+\xi)-f(t,\xi) \equiv g(t,z) \,,
\eeq
and the equation (\ref{55}) has three positive $p$-periodic solutions. However, $g(t,0)=0$ and $g_{zzz}(t,z)>0$ for $z>0$, contradicting the Lemma \ref{lma:5.3}. 
\medskip

Equations of the type (\ref{51}) occur often in ecological problems, see e.g., S. Ahmad and A.C. Lazer \cite{AS}, or P. Korman \cite{K}.


\begin{thebibliography}{99}
\bibitem{AS}
S. Ahmad and A.C. Lazer, Separated solutions of logistic equation with nonperiodic harvesting, {\em J. Math. Anal. Appl.} {\bf  445}, no. 1, 710-718  (2017).
\vspace{-0.2cm}

\bibitem{K}
P. Korman, Lectures on Differential Equations, AMS/MAA Textbooks, Volume $54$, 2019.
\vspace{-0.2cm}

\bibitem{KO} 
P. Korman and T. Ouyang,  Exact multiplicity results for two classes of periodic equations, {\em J. Math. Anal. Appl.} {\bf 194}, no. 3, 763-779  (1995).
\vspace{-0.2cm} 

\bibitem{P1}
A.A. Panov,  On the number of periodic solutions of polynomial differential equations. (Russian) {\em  Mat. Zametki} {\bf 64}, no. 5, 720-727  (1998); translation in {\em Math. Notes}, no. 5-6, 622-628 (1999).
\vspace{-0.2cm}

\bibitem{P}
V.A. Pliss,  Nonlocal Problems of the Theory of Oscillations. Translated from the Russian by Scripta Technica, Inc. Academic Press, New York-London (1966). 
\vspace{-0.2cm}

\bibitem{S}
A.  Sandqvist and K.M. Andersen,  On the number of closed solutions to an equation $x'=f(t,x)$, where $f_{x^n}(t,x)\geq 0$ ($n=1$,$2$, or $3$), {\em J. Math. Anal. Appl.} {\bf  159}, no. 1, 127-146  (1991).

\end{thebibliography}
\end{document}